\input amstex
\documentstyle{amsppt}
\magnification=\magstep1
 \hsize 13cm \vsize 18.35cm \pageno=1
\loadbold \loadmsam
    \loadmsbm
    \UseAMSsymbols
\topmatter
\NoRunningHeads
\title Note on the modified $q$-Bernstein polynomials
\endtitle
\author
  Taekyun Kim, Lee-Chae Jang, and Heungsu Yi
\endauthor
 \keywords : Bernstein polynomial, Bernstein operator,
  Bernoulli polynomial, generating function, Laurent series, shift difference operator
\endkeywords

\abstract  In the present paper, we propose the modified $q$-Bernstein polynomials of
degree $n$, which are different $q$-Bernstein polynomials of Phillips(see [4]). From these the
modified $q$-Bernstein polynomials of degree $n$, we derive some interesting recurrence formulae
for the modified $q$-Bernstein polynomials.
\endabstract
\thanks  2000 AMS Subject Classification: 11B68, 11S80
\newline
\endthanks
\endtopmatter

\document

{\bf\centerline {\S 1. Introduction}}

 \vskip 15pt

 Let $C[0,1]$ denote the set of continuous function on $[0,1]$. In [7], Bernstein introduced the
 following well-known linear positive operators.
 $$\eqalignno{ B_n (f:x)  &= \sum_{k=0}^n f\left( \frac{k}{n}\right) \binom{n}{k} x^k (1-x)^{n-k}\cr
 &=\sum_{k=0}^n f\left( \frac{k}{n}\right) B_{k,n}(x),  & (1) \cr}$$
 for $f\in C[0,1]$.
 $B_n(f:x)$ is called the Bernstein operator for $f$. The Bernstein polynomial of degree $n$ is defined by
 $$\eqalignno{ &B_{k,n} (x)  = \binom{n}{k} x^k (1-x)^{n-k},  & (2) \cr}$$
 for $k,n \in \Bbb Z_+$, where $x\in [0,1]$ and $\binom{n}{k}= \frac{n(n-1)\cdots (n-k+1)}{k!}$.
  It is easy to show that
 $$ B_{0,1}(x)=1-x, B_{1,1}(x)=x,$$
 $$ B_{0,2}(x)=(1-x)^2, B_{1,2}(x)= 2x(1-x), B_{2,2}(x)=x^2,$$
 $$B_{0,3}(x)=(1-x)^3, B_{1,3}(x)=3x(1-x)^2, B_{2,3}(x)=3x^2(1-x),$$
 $$ B_{3,3}(x)=x^3, \cdots.$$

Many researchers have studied the Bernstein polynomials in the area of approximation theory(see[1-8]).
For $k\in \Bbb Z_+$, it is easy to show that
 $$\eqalignno{ \frac{t^k e^{(1-x)t}x^k}{k!} &= \frac{x^k}{k!}\left( t^k\sum_{n=0}^\infty \frac{(1-x)^nt^n}{n!}\right)
 =\frac{x^k}{k!}\sum_{n=0}^\infty \frac{(1-x)^nt^{n+k}}{n!} \cr
 &= \frac{x^k}{k!}\sum_{n=0}^\infty \frac{(1-x)^n(n+1)\cdots (n+k)}{(n+k)!}t^{n+k}
 =x^k \sum_{n=k}^\infty \frac{(1-x)^{n-k} \binom{n}{k} t^n}{n!} \cr
 &=\sum_{n=k}\left( x^k (1-x)^{n-k}\binom{n}{k} \right) \frac{t^n}{n!}=\sum_{n=k}^\infty B_{k,n}(x) \frac{t^n}{n!},  &  \cr}$$
 and $B_{k,0}(x)=B_{k,1}(x)=\cdots =B_{k,k-1}(x)=0$. Thus, we obtain the generating function for $B_{k,n}(x)$ as follows:
 $$\eqalignno{ &F^{(k)} (t,x)= \frac{x^k e^{(1-x)t}t^k}{k!} =\sum_{n=0}^\infty B_{k,n}(x) \frac{t^n}{n!},  & (3) \cr}$$
 where $k\in \Bbb Z_+$ and $x\in [0,1]$. From (3), we can derive
$$ \eqalignno{& B_{k,n}(x) =  \cases \binom{n}{k}x^k (1-x)^{n-k}  & \text{if } n\geq k \\
0 & \text{if } n<k,  \endcases  & \cr}$$
for $n,k\in \Bbb Z_+$.

Let $q$ be regarded as a real number with $0<q<1$ and let us define the $q$-number as follows:
$$ [x]_q=[x:q]=\frac{1-q^x}{1-q}, \;\;\;\;\; \text{(see [1-7] )} .$$
Note that $\lim_{q\rightarrow 1}[x]_q=x$. In [4], Phillips introduced the $q$-extension
of Bernstein polynomials. Recently, Simsek and Acikgoz have also stdied the $q$-extension
of Bernstein type polynomials([5]). Their $q$-Bernstein type polynomials are given by
$$\eqalignno{ Y_n (k;x:q)& =\binom{n}{k} \frac{(-1)^k k!}{(1-q)^{n-k}} \sum_{m,l=0}^\infty \sum_{j=0}^{n-k}\binom{k+l-1}{l}\binom{n-k}{k}\cr &\;\;\;\;\;\;\; \times\left( \frac{(-1)^j q^{l+j(1-x)}S(m,k)(x\ln q)^m}{m!} \right),& \cr}$$
where $S(m,k)$ are the second kind stirling number.

In this paper we consider the $q$-extension of the generating function of Bernstein polynomials(see Eq.(3)).
From these $q$-extension of generating function for the Bernstein polynomials, we propose
the modified $q$-Bernstein polynomials of degree $n$, which are different $q$-Bernstein polynomials of Phillips.
By using the properties of the modified $q$-Bernstein polynomials, we can obtain some interesting recurrence formulae
for the modified $q$-Bernstein polynomials of degree $n$.

\vskip 10pt

{\bf\centerline {\S 2. The modified $q$-Bernstein polynomials}}

\vskip 10pt

For $q\in \Bbb R$ with $0<q<1$, let us consider the $q$-extension of Eq.(3) as follows:
 $$\eqalignno{ F_q^{(k)} (t,x) &= \frac{t^k e^{[1-x]_qt} [x]_q^k}{k!} = \frac{[x]_q^k}{k!}\sum_{n=0}^\infty \frac{[1-x]_q^n}{n!} t^{n+k} \cr
 &=[x]_q^k \sum_{n=0}^\infty \binom{n+k}{k} [1-x]_q^n  \frac{t^{n+k}}{(n+k)!} \cr
 &=\sum_{n=k}^\infty \binom{n}{k} [x]_q^k [1-x]_q^{n-k}  \frac{t^{n}}{n!} ,  & (4) \cr}$$
where $k,n \in \Bbb Z_+$ and $x \in [0,1]$. Note that $ \lim_{q\rightarrow1}F_q^{(k)} (t,x)=F(t,x)$. By (4), we can define the modified $q$-Bernstein polynomials as follows:
$$\eqalignno{& F_q^{(k)} (t,x) = \frac{t^k e^{[1-x]_qt} [x]_q^k}{k!} =
 =\sum_{n=k}^\infty B_{k,n}(x,q) \frac{t^{n}}{n!} ,  & (5) \cr}$$
where  $k,n \in \Bbb Z_+$ and $x \in [0,1]$.
By comparing the coefficients on the both sides of (4) and (5), we obtain the following theorem.

\vskip 10pt

\proclaim {Theorem 1} For $k,n \in\Bbb Z_+, x \in [0,1]$, we have
$$ \eqalignno{& B_{k,n}(x,q) =  \cases \binom{n}{k}[x]_q^k [1-x]_q^{n-k}  & \text{if } n\geq k \\
0 & \text{if } n<k.  \endcases  & \cr}$$
\endproclaim

\vskip 10pt

For $0\leq k\leq n$, we have
$$\eqalignno{& [1-x]_q B_{k,n-1}(x,q)+[x]_qB_{k-1,n-1}(x,q)\cr
&= [1-x]_q\binom{n-1}{k} [x]_q^k [1-x]^{n-1-k}+[x]_q\binom{n-1}{k-1}[x]_q^{k-1}[1-x]_q^{n-k}\cr
&= \binom{n-1}{k} [x]_q^k [1-x]^{n-k}+\binom{n-1}{k-1}[x]_q^{k}[1-x]_q^{n-k}\cr
&=\binom{n}{k} [x]_q^k [1-x]_q^{n-k},  &  \cr}$$
and the derivative of the modified $q$-Bernstein polynomials of degree $n$ are also polynomials of degree $n-1$. That is,
$$\eqalignno{&\frac{d}{dx}B_{k,n}(x,q) \cr &= \binom{n}{k} k [x]_q^{k-1}[1-x]_q^{n-k} \frac{\ln q}{q-1}q^x
+\binom{n}{k} [x]_q^k (n-k)[1-x]_q^{n-k-1}\left(-\frac{\ln q}{q-1}\right) q^{1-x}\cr
&= \frac{\ln q}{q-1} \left\{ \binom{n}{k} k [x]_q^{k-1}[1-x]_q^{n-k}q^x-
\binom{n}{k} [x]_q^k (n-k)[1-x]_q^{n-k-1}q^{1-x}\right\}\cr
&= n \left( q^x B_{k-1,n-1}(x,q)-q^{1-x}B_{k,n-1}(x,q)\right)\frac{\ln q}{q-1}.  &  \cr}$$
Therefore, we obtain the following recurrence formulae.

\vskip 10pt

\proclaim {Theorem 2 (Recurrence formulae for $B_{k,n}(x,q)$)} For $k,n \in\Bbb Z_+, x \in [0,1]$, we have
$$  [1-x]_qB_{k,n-1}(x,q)+[x]_qB_{k-1,n-1}(x,q)= \binom{n}{k}[x]_q^k[1-x]_q^{n-k}=B_{k,n}(x,q), $$
and
$$\frac{d}{dx}B_{k,n}(x,q) =n \left(q^xB_{k-1,n-1}(x,q)-q^{1-x}B_{k,n-1}(x,q)\right)\frac{\ln q}{q-1}.$$
\endproclaim

\vskip 10pt

Let $f$ be a continuous function on $[0,1]$. Then the modified $q$-Bernstein operator is defined by
$$ \eqalignno{& B_{n,q}(f:x) = \sum_{j=0}^n f\left(\frac{j}{n}\right)  B_{j,n}(x,q),  & (6) \cr}$$
where $0\leq x\leq 1, n \in \Bbb Z_+$. By Theorem 1 and (6), we see that
$$ \eqalignno{ B_{n,q}(f:x) &=B_{n,q}f(x) = \sum_{k=0}^n f\left(\frac{k}{n}\right) \binom{n}{k} [x]_q^k[1-x]_q^{n-k} \cr
&=[x]_q (1-[1-x]_q [x]_q (q-1))^{n-1},  &  \cr}$$
where $f(x)=x$. Thus, we have
$$ \eqalignno{& B_{n,q}(f:x) = f([x]_q) \left( 1+(1-q)[x]_q[1-x]_q \right)^{n-1}.  & (7) \cr}$$
From Theorem 1, we note that
$$ \eqalignno{ \sum_{k=0}^n B_{k,n}(x,q) &=\sum_{k=0}^n \binom{n}{k} [x]_q^k[1-x]_q^{n-k} \cr
&=\sum_{k=0}^n \binom{n}{k} [x]_q^k (1-q^{1-x}[x]_q)^{n-k}\cr
&=(1+[x]_q [1-x]_q (q-1))^{n}=B_{n,q}(1:x).  &  \cr}$$

The modified $q$-Bernstein polynomials are symmetric polynomials. That is, by the definition of the modified $q$-Bernstein polynomials of degree $n$, we see that
$$B_{k,n}(x,q)=\binom{n}{k} [x]_q^k [1-x]_q^{n-k}.$$
Thus, we have
$$ B_{n-k,n}(1-x,q)=\binom{n}{n-k}[1-x]_q^{n-k}[x]_q^k =\binom{n}{k}[x]_q^k[1-x]_q^{n-k}.$$
Therefore, we obtain the following theorem.

\vskip 10pt

\proclaim {Theorem 3 } For $k,n \in\Bbb Z_+, x \in [0,1]$, we have
$$ B_{n-k,n}(1-x,q)=B_{k,n}(x,q), $$
and
$$ \sum_{k=0}^n B_{k,n}(x,q)=(1+[x]_q[1-x]_q(1-q))^n =B_{n,q}(1:x).$$
\endproclaim

\vskip 10pt

For $t \in \Bbb C, x\in [0,1]$, and $n \in \Bbb Z_+$, we consider
$$\eqalignno{& \frac{n!}{2\pi i} \oint_C \frac{([x]_qt)^k}{k!}e^{([1-x]_qt)} \frac{dt}{t^{n+1}}
 ,  & (9) \cr}$$
where $C$ is a circle around the origin and integration is in the positive direction.
By the definition of the modified $q$-Bernstein polynomials and Laurent series, we see that
$$\eqalignno{&  \oint_C \frac{([x]_qt)^k}{k!}e^{[1-x]_qt} \frac{dt}{t^{n+1}} \cr
&=\sum_{m=0}^\infty \oint_C \frac{B_{k,n}(x,q)t^m}{m!}\frac{dt}{t^{n+1}}=\frac{B_{k,n}(x,q)}{n!}2\pi i
 .  & (10) \cr}$$
From (9) and (10), we note that
$$\eqalignno{& \frac{n!}{2\pi i} \oint_C \frac{([x]_qt)^k}{k!}e^{[1-x]_qt} \frac{dt}{t^{n+1}}= B_{k,n}(x,q) .  & (11) \cr}$$
Also, we see that
$$\eqalignno{&  \oint_C \frac{([x]_qt)^k}{k!}e^{[1-x]_qt} \frac{dt}{t^{n+1}} \cr
&= \frac{[x]_q^k}{k!}  \sum_{m=0}^\infty \oint_C t^{m-n-1+k} dt \frac{[1-x]_q^m }{m!} \cr
&= \frac{[x]_q^k}{k!} \left( \frac{[1-x]_q^{n-k}}{(n-k)!}\right) 2\pi i\cr
&=\frac{[x]_q^k [1-x]_q^{n-k}}{k!(n-k)!} 2\pi i
 .  & (12) \cr}$$
 From (9) and (12), we have
 $$\eqalignno{&  \frac{n!}{2\pi i} \oint_C \frac{([x]_qt)^k}{k!}e^{([1-x]_qt)} \frac{dt}{t^{n+1}} 
= \binom{n}{k} [x]_q^k [1-x]_q^{n-k}. & (13) \cr}$$
By (11) and (13), we easily see that
 $$\eqalignno{&  B_{k,n}(x,q)=\binom{n}{k} [x]_q^k [1-x]_q^{n-k} . & (14) \cr}$$
From (14), we derive
 $$\eqalignno{& \left(  \frac{n-k}{n} \right) B_{k,n}(x,q) + \left(  \frac{k+1}{n} \right) B_{k+1,n}(x,q) \cr
 &= \frac{n-k}{n} \binom{n}{k} [x]_q^k [1-x]_q^{n-k} + \frac{k+1}{n} \binom{n}{k+1} [x]_q^{k+1}[1-x]_q^{n-k-1} \cr
 &= \frac{(n-1)!}{k!(n-k-1)!}[x]_q^k[1-x]_q^{n-k}+\frac{(n-1)!}{k!(n-k-1)!}[x]^{k+1}[1-x]_q^{n-k-1} \cr
 &= [1-x]_qB_{k,n-1}(x,q)+[x]_qB_{k,n-1}(x,q)\cr
 &=B_{k,n-1}(x,q)+ [x]_q(1-q^{1-x})B_{k,n-1}(x,q) \cr
 &=B_{k,n-1}(x,q) +(1-q)[x]_q[1-x]_qB_{k,n-1}(x,q). &  \cr}$$
 Therefore, we can write the modified $q$-Bernstein polynomials as a linear combination of polynomials of higher order as follows:

\vskip 10pt

\proclaim {Theorem 4 } For $k \in\Bbb Z_+, n \in \Bbb N$, and $x \in [0,1]$, we have
$$ \eqalignno{&\left( \frac{n-k}{n}\right)B_{k,n}(x,q)+  \left( \frac{k+1}{n}\right)B_{k+1,n}(x,q)\cr
&= B_{k, n-1}(x,q) + (1-q)[x]_q[1-x]_q B_{k,n-1}(x,q). &  \cr}$$
\endproclaim

\vskip 10pt

By (14), we easily see that
 $$\eqalignno{& \left( \frac{n-k+1}{k} \right) \left( \frac{[x]_q}{[1-x]_q} \right) B_{k-1,n}(x,q) \cr
 &=  \left( \frac{n-k+1}{k} \right) \left( \frac{[x]_q}{[1-x]_q} \right) \binom{n}{k-1} [x]_q^{k-1}[1-x]_q^{n_k+1} \cr
 &= \frac{n!}{k!(n-k)!} [x]_q^k [1-x]_q^{n-k} =\binom{n}{k} [x]_q^k [1-x]_q^{n-k}. & \cr}$$
Thus, we obtain the following corollary.

\vskip 10pt

\proclaim {Corollary 5 } For $n, k \in \Bbb N$, and $x \in [0,1]$, we have
$$ \left( \frac{n-k+1}{k}\right) \left( \frac{[x]_q}{[1-x]_q}\right)B_{k-1,n}(x,q)
= B_{k, n}(x,q) . $$
\endproclaim

\vskip 10pt

From the definition of the modified $q$-Bernstein polynomials and binomial theorem, we note that 
 $$\eqalignno{B_{k,n}(x,q)&= \binom{n}{k} [x]_q^k [1-x]_q^{n-k}=\binom{n}{k}[x]_q^k (1-q^{1-x}[x]_q)^{n-k} \cr
 &=\binom{n}{k} [x]_q^k \sum_{l=0}^{n-k} \binom{n-k}{l} (-1)^l q^{l(1-x)} [x]_q^l\cr
 &= \sum_{l=0}^{n-k}\binom{n}{k} \binom{n-k}{l} (-1)^l q^{l(1-x)} [x]_q^{l+k} \cr
 &= \sum_{l=0}^{n-k}\binom{k+l}{k} \binom{n}{k+l} (-1)^l q^{l(1-x)} [x]_q^{l+k} \cr
 &= \sum_{i=k}^{n}\binom{i}{k} \binom{n}{i} (-1)^{i-k} q^{(1-x)(i-k)} [x]_q^{i}. &  \cr}$$
Therefore, we obtain the following theorem.

\vskip 10pt

\proclaim {Theorem 6 } For $k, n \in \Bbb Z_+$, and $x \in [0,1]$, we have
$$ B_{k,n}(x,q)=\sum_{i=k}^{n}\binom{i}{k} \binom{n}{i} (-1)^{i-k} q^{(1-x)(i-k)} [x]_q^{i}.  $$
\endproclaim

\vskip 10pt

It is possible to write each power basis element of $[x]_q^k$, in the linear combination of the modified $q$-Bernstein polynomials by using the degree evaluation formulae and induction method in mathematics. From the property of the modified $q$-Bernstein polynomials, we easily see that
 $$\eqalignno{ \sum_{k=0}^n \frac{k}{n} B_{k,n}(x,q) &= \sum_{k=0}^{n} \binom{n}{k} [x]_q^k [1-x]_q^{n-k}\cr
 &= \sum_{k=1}^{n} \binom{n-1}{k-1} [x]_q^k [1-x]_q^{n-k}\cr
 &= \sum_{k=0}^{n-1} \binom{n-1}{k} [x]_q^{k+1} [1-x]_q^{n-1-k}\cr
 &= [x]_q ([x]_q+[1-x]_q)^{n-1}, &  \cr}$$
 and
  $$\eqalignno{ \sum_{k=1}^n \frac{\binom{k}{2}}{\binom{n}{2}} B_{k,n}(x,q) 
  &= \sum_{k=1}^{n} \frac{k(k-1)}{n(n-1)}\binom{n}{k} [1-x]_q^{n-k}[x]_q^k \cr
 &= \sum_{k=2}^\infty \frac{k(k-1)}{n(n-1)} \binom{n}{k} [x]_q^k [1-x]_q^{n-k}\cr
 &= \sum_{k=2}^n \binom{n-2}{k-2} [x]_q^{k} [1-x]_q^{n-k}\cr
 &=\sum_{k=0}^{n-2} \binom{n-2}{k} [x]_q^{k+2} [1-x]_q^{n-2-k}\cr
 &= [x]_q^2 ([x]_q+[1-x]_q)^{n-2}. &  \cr}$$
Continuing this process, we obtain 
 $$\eqalignno{ \sum_{k=i-1}^n \frac{\binom{k}{i}}{\binom{n}{i}} B_{k,n}(x,q)
  = [x]_q^i ([x]_q+[1-x]_q)^{n-i}, &  \cr}$$
 for $i \in \Bbb N$.
 Therefore we obtain the following theorem.

\vskip 10pt

\proclaim {Theorem 7 } For $n \in \Bbb Z_+$, $i \in \Bbb N$ and $x \in [0,1]$, we have
$$ \frac{1}{([1-x]_q+[x]_q)^{n-i}}\sum_{k=i-1}^n \frac{\binom{k}{i}}{\binom{n}{i}} B_{k,n}(x,q)
  = [x]_q^i  .  $$
\endproclaim

\vskip 10pt
 
The Bernoulli polynomials of order $k (\in \Bbb N)$ are defined as
$$ \eqalignno{&
\left( \frac{t}{e^t -1} \right)^k e^{xt} =
\underbrace{\left(\frac{t}{e^t -1} \right)^k \times \cdots \times \left(\frac{t}{e^t -1} \right)^k }_{k-\text{times}}
e^{xt} = \sum_{n=0}^\infty B_n^{(k)} (x) \frac{t^n}{n!},  & (15) \cr}
$$
and $ B_n^{(k)}=B_n^{(k)}(0)$ are called the $n$-th Bernoulli numbers of order $k$. 
It is well known that the second kind stirling numbers are defined as
 $$\eqalignno{&\frac{(e^t-1)^k}{k!} =\sum_{n=0}^\infty S(n,k) \frac{t^n}{n!},  & (16) \cr}$$
for $k \in \Bbb N$. From (5), we note that
$$\eqalignno{& \frac{([x]_qt)^k e^{[1-x]_qt}}{k!}=\frac{[x]_q^k(e^t-1)^k}{k!}\left( \frac{t}{e^t-1}\right)^k e^{[1-x]_qt} \cr
&=[x]_q^k \left( \sum_{m=0}^\infty S(m,k) \frac{t^m}{m!}\right) \left(\sum_{n=0}^\infty B_n^{(k)} ([1-x]_q)\frac{t^n}{n!}\right) \cr
&= [x]^k \sum_{l=0, m+n=l}^\infty \left( \sum_{n=0}^l \frac{B_n^{(k)}([1-x]_q)S(l-n,k)\binom{l}{n}}{n!(l-n)!}\right) \frac{t^l}{l!}
. & (17) \cr}$$
By (5) and (17), we have
$$B_{k,l}(x,q)=[x]_q^k \sum_{n=0}^l B_n^{(k)} ([1-x]_q)S(l-n,k)\binom{l}{n},$$
and $B_{k,0}(x,q)=B_{k,1}(x,q)= \cdots =B_{k,k-1}(x,q)=0$, where $B_n^{(k)}([1-x]_q)$ are called the $n$-th Bernoulli polynomials of order $k$.

Let $\Delta$ be the shift difference operator with $\Delta f(x)=f(x+1)-f(x)$. By iterative method, we easily see that
$$ \eqalignno{& \Delta^n f(0)= \sum_{k=0}^n \binom{n}{k} (-1)^{n-k} f(k),  & (18) \cr}$$ 
for $n\in \Bbb N$. From (16) and (18), we note that
$$ \eqalignno{  \sum_{n=0}^\infty S(n,k) \frac{t^n}{n!}& =\frac{1}{k!}\sum_{l=0}^k \binom{k}{l} (-1)^{k-l}e^{lt}\cr
&=\sum_{n=0}^\infty \left\{  \frac{1}{k!} \sum_{l=0}^k \binom{k}{l} (-1)^{k-l} l^n  \right\} \frac{t^n}{n!} \cr
&= \sum_{n=0}^\infty \frac{\Delta^k 0^n}{k!}\frac{t^n}{n!}.
 &  \cr}$$ 
By comparing the coefficients on the both sides, we have
$$ \eqalignno{&   S(n,k)= \frac{\Delta^k 0^n}{k!}, & (19) \cr}$$ 
for $n,k\in \Bbb Z_+$. Thus, we note that
$$ \eqalignno{& B_{k,l}(x,q)=[x]_q^k \sum_{n=0}^l B_n^{(k)} ([1-x]_q)\binom{l}{n} \frac{\Delta^k 0^{l-n}}{k!}.  & (20) \cr}$$
Let $(Eh)(x)=h(x+1)$ be the shift operator. Then the $q$-difference operator is defined by
$$ \eqalignno{& \Delta_q^n = \Pi_{i=0}^{n-1} (E-q^i I), \text{ (see [2])}, &  \cr}$$
where $I$ is an identity operator. For $f\in C[0,1]$ and $n\in \Bbb N$, we have 
$$ \eqalignno{& \Delta_q^n f(0)= \sum_{k=0}^n \binom{n}{k}_q (-1)^k q^{\binom{n}{2}} f(n-k),  &  \cr}$$
where $\binom{n}{k}_q$ is Gaussian binomial coefficient. 

Let $F_q(t)$ be the generating function of the $q$-extension of the second kind stirling number as follows:
$$ \eqalignno{F_q(t) &=\frac{q^{-\binom{k}{2}}}{[k]_q!} \sum_{0\leq j \leq k} (-1)^{k-j}\binom{k}{j}_q q^{\binom{k-j}{2}} e^{[j]_qt} \cr
&=\sum_{n=0}^\infty S(n,k:q) \frac{t^n}{n!}, \text{ (see [2])}.  & (21) \cr}$$
From (21), we have
$$ \eqalignno{S(n,k:q) &= \frac{q^{-\binom{k}{2}}}{[k]_q!} \sum_{j=0}^k (-1)^j q^{\binom{j}{2}} \binom{k}{j}_q [k-j]_q^n \cr
 &= \frac{q^{-\binom{k}{2}}}{[k]_q!} \Delta_q^k 0^n, & (22)  \cr}$$
where $[k]_q!=[k]_q [k-1]_q \cdots [2]_q [1]_q$. It is not difficult to show that
$$ \eqalignno{& [x]_q^n = \sum_{k=0}^n q^{\binom{k}{2}} \binom{x}{k}_q [k]_q! S(n,k:q), \text{ (see [2])}, & (23)  \cr}$$
Thus, we have 
$$ \eqalignno{& \sum_{k=0}^i q^{\binom{k}{2}} \binom{x}{k}_q [k]_q! S(i,k:q)= 
\frac{1}{([1-x]_q+[x]_q)^{n-i}}\sum_{k=i-1}^n \frac{\binom{k}{i}}{\binom{n}{i}} B_{k,n}(x,q). & (24) \cr}$$ 
Therefore, we obtain the following theorem.

\vskip 10pt

\proclaim {Theorem 8 } For $n \in \Bbb Z_+$, $i \in \Bbb N$ and $x \in [0,1]$, we have
$$ \frac{1}{([1-x]_q+[x]_q)^{n-i}}\sum_{k=i-1}^n \frac{\binom{k}{i}}{\binom{n}{i}} B_{k,n}(x,q)
  = \sum_{k=0}^i q^{\binom{k}{2}} \binom{x}{k}_q [k]_q! S(i,k:q),  $$
where $\binom{x}{k}_q=\frac{[x]_q[x-1]_q \cdots [x-k+1]_q}{[k]_q!}$.
\endproclaim

\vskip 20pt

\quad Taekyun Kim

\quad  Division of General Education,

\quad Kwangwoon  University,

\quad Seoul 139-701, South Korea
 \quad e-mail:\text{tkkim\@kw.ac.kr}

\vskip 10pt

\quad Lee-Chae Jang

\quad Department of Mathematics and Comp. Sci.,

 \quad  Kon-Kuk University,

 \quad Chungju 138-701, Korea
 \quad e-mail:\text{leechae.jang\@kku.ac.kr}

\vskip 10pt
\quad Heungsu Yi

\quad  Department of Mathematic,

\quad Kwangwoon  University,

\quad Seoul 139-701, South Korea
 \quad e-mail:\text{hsyi\@kw.ac.kr}

\enddocument